%% file: r-dim_lambda_max_fin_v6.tex
\documentclass[a4paper,reqno,12pt]{amsart}
\usepackage{rotating,pdflscape,color,amssymb,hyperref,subfigure,psfrag,amsmath,eufrak,bbm,cancel}
\author[Florent Benaych-Georges]{Florent Benaych-Georges}\address{Florent Benaych-Georges, LPMA,  UPMC Univ Paris 6, Case courier 188, 4, Place Jussieu, 75252 Paris Cedex 05, France, and CMAP, \'Ecole Polytechnique, route de Saclay, 91128 Palaiseau Cedex, France.} \email{florent.benaych@upmc.fr}
\urladdr{http://www.cmapx.polytechnique.fr/\~\/benaych/}

\author[Raj Rao Nadakuditi]{Raj Rao Nadakuditi}\address{Raj Rao Nadakuditi, Department of Electrical Engineering and Computer Science, University of Michigan, 1301 Beal Avenue, Ann Arbor, MI 48109. USA.}
\email{rajnrao@eecs.umich.edu}
\urladdr{http://www.eecs.umich.edu/\~\/rajnrao/}

\title[Low rank perturbations of large random matrices]{The eigenvalues and eigenvectors of finite, low rank perturbations of large random matrices}

\thanks{F.B.G's work was partially supported by the \emph{Agence Nationale de la Recherche} grant ANR-08-BLAN-0311-03. R.R.N's research was partially supported by an Office of Naval Research postdoctoral fellowship award and grant N00014-07-1-0269. R.R.N thanks Arthur Baggeroer for his feedback, support and encouragement. We thank Alan Edelman for feedback and encouragement and for facilitating this collaboration by hosting F.B.G's stay at M.I.T. We gratefully acknowledge the Singapore-MIT alliance for funding F.B.G's stay.}

\date{\today}
\newcommand{\tta}{\theta}

\newcommand{\ovl}{\overline}
\newcommand{\bbm}{\begin{bmatrix}}
\newcommand{\ebm}{\end{bmatrix}}
\newcommand{\bes}{\begin{equation*}}
\newcommand{\ees}{\end{equation*}}
\newcommand{\be}{\begin{equation}}
\newcommand{\ee}{\end{equation}}
\newcommand{\beqy}{\begin{eqnarray}}
\newcommand{\eeqy}{\end{eqnarray}}
\newcommand{\beq}{\begin{eqnarray*}}
\newcommand{\eeq}{\end{eqnarray*}}

\newcommand{\lan}{\langle}
\newcommand{\ran}{\rangle}

\newcommand{\diag}{\operatorname{diag}}
\newcommand{\Diag}{\operatorname{diag}}

\newcommand{\Pro}{\mathbb{P}}

\newcommand{\supp}{\operatorname{supp}}

\newcommand{\Tr}{\operatorname{Tr}}

\newcommand{\ninf}{\underset{n\to\infty}{\longrightarrow}}

\newcommand{\one}{\mathbbm{1}}

\newcommand{\E}{\mathbb{E}}

\newcommand{\R}{\mathbb{R}}
\newcommand{\C}{\mathbb{C}}

\newcommand{\ud}{\mathrm{d}}

\newcommand{\pro}{probability }
\newcommand{\ie}{\textit{i.e.}}
\newcommand{\f}{\frac}
\newcommand{\ff}{\frac{1}}
\newcommand{\lf}{\left}
\newcommand{\ri}{\right}

\newcommand{\st}{such that }
\newcommand{\la}{\lambda}

\newcommand{\ste}{\, ;\, }

\newcommand{\eps}{\varepsilon}

\newcommand{\bxp}{\boxplus}
\newcommand{\bxt}{\boxtimes}

\newcommand{\wtX}{\widetilde{X}_n}
\newcommand{\wtl}{\widetilde{\la}}
\newcommand{\convas}{\overset{\textrm{a.s.}}{\longrightarrow}}

\newcommand{\lto}{\longrightarrow}
\newcommand{\al}{\alpha}
\newcommand{\bck}{\backslash}

\newtheorem{Th}{Theorem}[section]

\newtheorem{propo}[Th]{Proposition}

\newtheorem{lem}[Th]{Lemma}

\newtheorem{rmq}[Th]{Remark}

\newenvironment{pr}{\noindent {\it Proof. }}{\hfill$\square$}

\long\def\symbolfootnote[#1]#2{\begingroup
\def\thefootnote{\fnsymbol{footnote}}\footnote[#1]{#2}\endgroup}

\begin{document}
\maketitle
\begin{abstract}
We consider the eigenvalues and eigenvectors of finite,
low rank perturbations of random matrices.
Specifically, we prove almost sure convergence of the extreme eigenvalues and appropriate
projections of the corresponding eigenvectors of the perturbed matrix for additive and multiplicative perturbation models.

The limiting non-random value is shown to depend explicitly on the limiting eigenvalue distribution of the unperturbed random matrix and the assumed perturbation model via integral
transforms that correspond to very well known objects in free
probability theory that linearize non-commutative free additive and
multiplicative convolution. Furthermore, we uncover a phase
transition phenomenon whereby the large matrix limit of the extreme
eigenvalues of the perturbed matrix differs from that of the original
matrix if and only if the eigenvalues of the perturbing matrix are
above a certain critical threshold. Square root decay of the eigenvalue density at the edge is sufficient to ensure that this threshold is finite.
This critical threshold is
intimately related to the same aforementioned integral transforms and our proof techniques bring this connection and the origin of the phase transition into   focus. Consequently, our results extend the class of `spiked' random matrix models about which such predictions (called the {\it BBP phase transition}) can be made well beyond the Wigner, Wishart and Jacobi random ensembles found in the literature.
We examine the impact of this eigenvalue phase transition on the associated eigenvectors and observe an analogous phase transition in the eigenvectors. Various extensions of our results to the problem of
non-extreme eigenvalues are discussed.
\end{abstract}


\input{introduction_fin_v6}
\input{mainresults_fin_v6}

\input{examples_fin_v6}

\input{sketch-proofs_fin_v6}
\input{masterequations_fin_v6}

\input{proofsquareadd_fin_v6}

\input{convergencefacts_fin_v6}


\end{document}

%% file: introduction_fin_v6.tex
  \section{Introduction}
Let $X_n$ be an $n \times n$ symmetric (or Hermitian)
 matrix with   eigenvalues $\la_1(X_n), \ldots, \la_n(X_n)$ and $P_n$ be an $n \times n$ symmetric (or Hermitian)
 matrix with rank $r \leq n$ and   non-zero eigenvalues $\theta_1, \ldots, \theta_r$. A fundamental question in matrix analysis is the following \cite{bns78,agg88}:
\begin{quote}
How are the eigenvalues and eigenvectors of $X_n + P_n$ related to the eigenvalues and eigenvectors of $X_n$ and $P_n$?
\end{quote}
When $X_n$ and $P_n$ are diagonalized by the same eigenvectors, we have $\la_i(X_n+P_n) = \la_j(X_n) + \la_k(P_n)$ for appropriate choice of indices $i,j,k \in \{1,\ldots,n\}$. In the general setting, however, the answer is complicated by the fact that the eigenvalues and eigenvectors of their sum depend on the relationship between the eigenspaces of the individual matrices.

In this scenario, one can use Weyl's interlacing inequalities and Horn inequalities \cite{hj85} to obtain coarse bounds for the eigenvalues  of the sum in terms of the eigenvalues of $X_n$. When the norm of $P_n$ is small relative to the norm of $X_n$, tools from perturbation theory (see \cite[Chapter 6]{hj85} or \cite{ss90}) can be employed to improve the characterization of the bounded set in which the eigenvalues of the sum must lie. Exploiting any special structure in the matrices allows us to refine these bounds \cite{in09} but this is pretty much as far as the theory goes. Instead of exact answers we must resort to a system of coupled  inequalities. Describing the behavior of the eigenvectors of the sum is even more complicated.

\begin{figure}
\vspace{0.50in}
\centering
\subfigure[Largest eigenvalue $\rho>b$ in blue when $\theta > \theta_c$]{
\psfrag{c}{$\substack{\\ \\ \rho}$}
\psfrag{a}{$\substack{\\ \\ a}$}
\psfrag{b}{$\substack{\\ \\ b}$}
\includegraphics[width=2.4in]{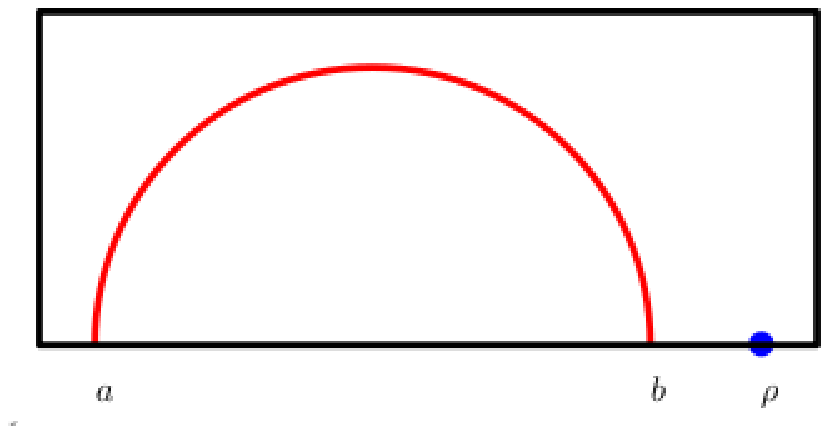}
\label{fig:evphaseon}
}\hspace{1.2in}
\subfigure[Associated eigenvector when $\theta > \theta_c$]{
\psfrag{u}{$u$}
\psfrag{u1}{$\widetilde{u}$}
\includegraphics[width=1.7in]{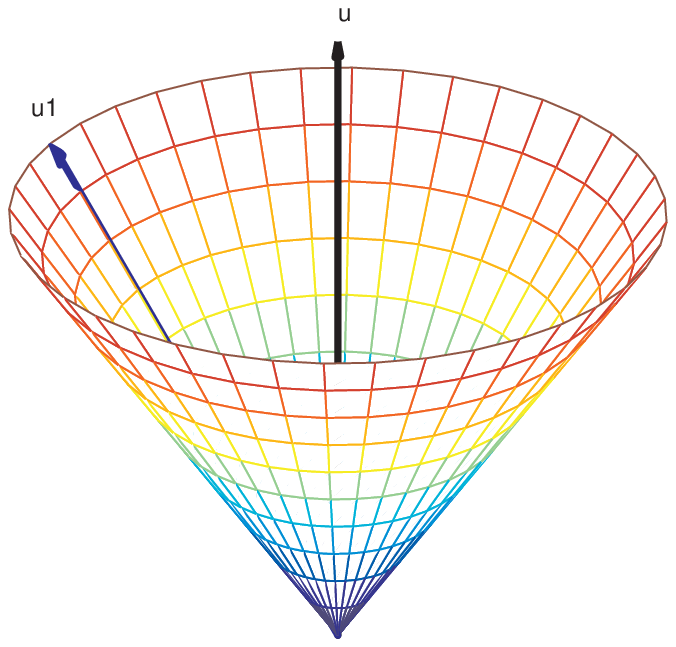}
\label{fig:eigproj}}
\\[0.5in]
\subfigure[Largest eigenvalue $\rho=b$ in blue when $\theta \leq \theta_c$]{
\psfrag{a}{$\substack{\\ \\ a}$}
\psfrag{b}{$\substack{\\ \\ \rho=b}$}
\includegraphics[width=2.4in]{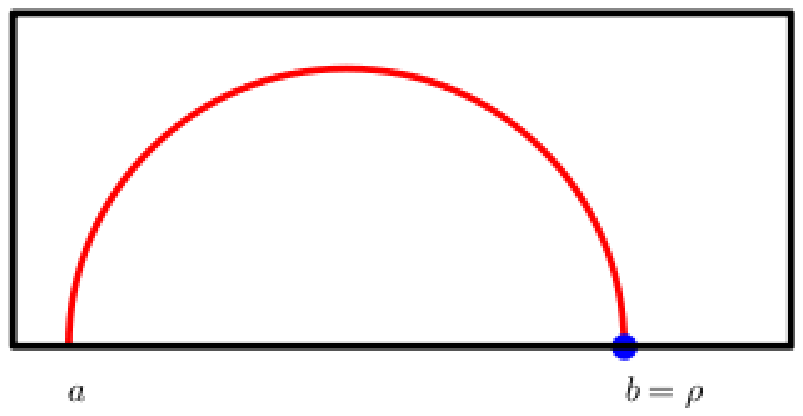}
\label{fig:evphaseoff}
}\hspace{0.8in}
\subfigure[Associated eigenvector when $\theta \leq \theta_c$]{
\psfrag{u}{$u$}
\psfrag{u1}{$\widetilde{u}$}
\includegraphics[width=1.7in]{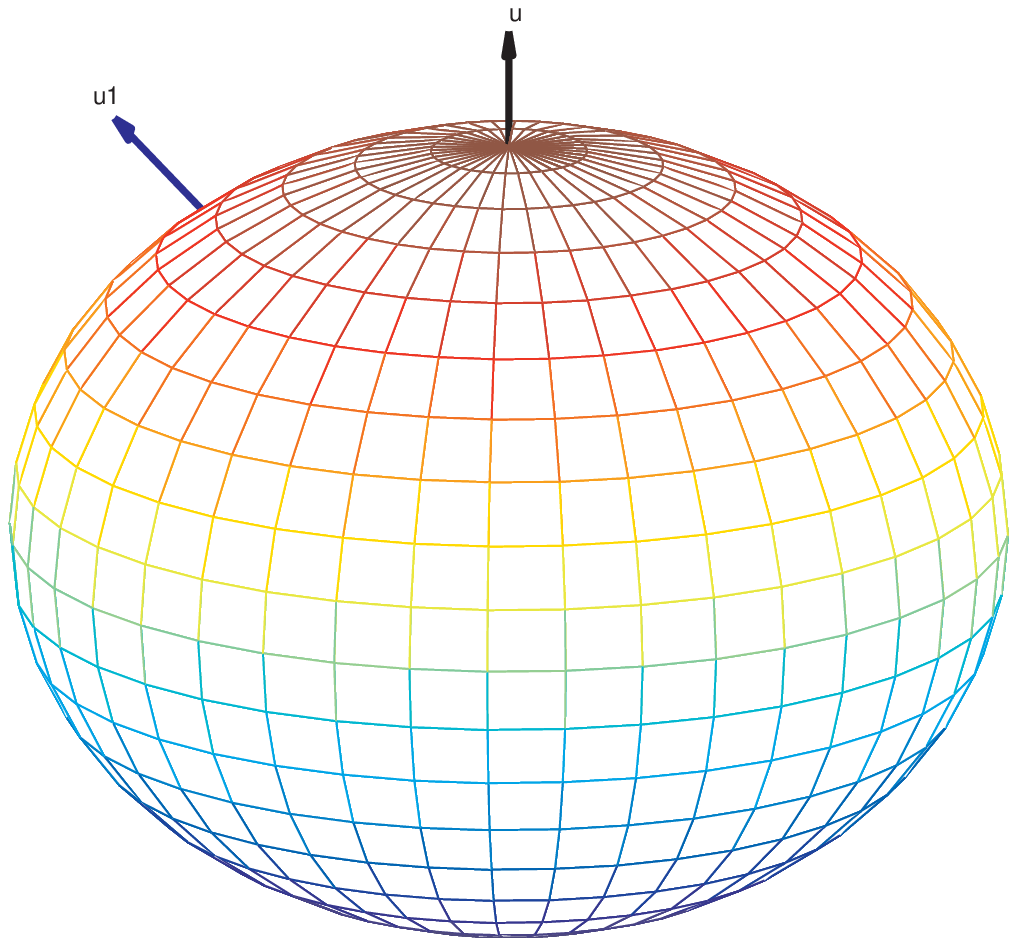}
\label{fig:eigsphere}
}
\caption{Assume that the limiting eigenvalue distribution of $X_n$ is $\mu_X$   with largest eigenvalue $b$. Consider the matrix $P_n:=\theta u u^*$ with rank $r = 1$ and largest eigenvalue $\theta$($>0$, say). The vector $u$ is an $n \times 1$ vector chosen uniformly at random from the unit $n$-sphere. The largest eigenvalue of $X_n + P_n$ will differ from $b$  if and only if $\theta$ is greater than some critical value $\theta_c$. In this event, the largest eigenvalue will be concentrated around $\rho$ with high probability as in  (a). The associated eigenvector $\widetilde{u}$ will, with high probability, lie on a cone around $u$ as in (b). When $\theta \leq \theta_c$, a phase transition occurs so that with high probability, the largest eigenvalue of the sum will equal $b$ as in (c) and the corresponding eigenvector will be uniformly distributed on the unit sphere as in (d).}
\label{fig:summary}
\end{figure}

Surprisingly, adding some randomness to the eigenspaces permits further analytical progress. Specifically, if the eigenspaces are assumed to be ``in generic position with respect to each other'', then in place of eigenvalue bounds we have simple, exact answers that are to be interpreted probabilistically. These results bring into   focus a   phase transition phenomenon of the kind illustrated in Figure \ref{fig:summary} for the eigenvalues and eigenvectors of $X_n + P_n$ and $X_n\times(I_n+P_n)$. A precise statement of the results may be found in Section \ref{sec:main results}.

Examining the structure of the analytical expression for the critical values $\theta_c$ and $\rho$ in Figure \ref{fig:summary} reveals a common underlying theme in the additive and multiplicative perturbation settings. The critical values $\theta_c$ and $\rho$ in Figure \ref{fig:summary} are related to integral transforms of the limiting eigenvalue distribution $\mu_X$ of $X_n$.  It turns out that the 
 integral transforms that emerge in the respective additive and multiplicative cases are deeply related to very well known objects in free probability theory \cite{vdn91,hp00,agz09} that linearize
free additive and multiplicative   convolutions respectively.
In a forthcoming paper \cite{flo-raj-rectangular}, we consider the analogue of the problem for the extreme singular values of finite rank deformations of rectangular random matrices. There too, a phase transition occurs at a threshold determined by an integral transform which plays an analogous role in the computation  of the rectangular additive free convolution \cite{b09,b09R=C,bg.sph.int}.
The emergence of these transforms in the context of the study of the extreme or isolated eigenvalue behavior should be of independent interest to free probabilists.
%
 In doing so, we extend the results found in the literature  about the  so-called {\it BBP phase transition} (named after Baik, Ben Arous, P\'ech\'e because of their seminal paper \cite{bbp05})  for the eigenvalue phase transition in such finite, low rank perturbation models well beyond the Wigner \cite{sp06,hr07,pf07,cdf09,bff09}, Wishart \cite{bbp05,bs06, e07,p07,n08} and Jacobi settings \cite{r09}. In our situation, the distribution $\mu_X$ in Figure \ref{fig:summary} can be \textit{any} probability measure. Consequently, the aforementioned results in the literature can be rederived rather simply using the formulas in Section \ref{sec:main results} by substituting $\mu_X$ with the semi-circle measure \cite{w58} (for Wigner matrices), the Mar\v{c}enko-Pastur measure \cite{mp67} (for Wishart matrices) or the free Jacobi measure (for Jacobi matrices \cite{c05}). Concrete computations are presented in Section \ref{sec:examples}.
 
The development of the eigenvector aspect   is another  contribution that we would like to highlight. Generally speaking, the eigenvector question has received   less   attention in random matrix theory and in free probability theory. A notable exception is the recent body of work on the eigenvectors of spiked Wishart matrices \cite{p07, hr07, n08} which corresponds to $\mu_X$ being the Mar\v{c}enko-Pastur measure. In this paper, we extend their results for multiplicative models of the kind $(I+P_n)^{1/2}X_n(I+P_n)^{1/2}$ to the setting where $\mu_X$ is an arbitrary probability measure and obtain new results for the eigenvectors for additive models of the form $X_n + P_n$. 

Our proofs rely  on the derivation of master equation representations of the eigenvalues and eigenvectors of the perturbed matrix and the subsequent application of concentration inequalities for random vectors uniformly distributed on high dimensional unit spheres (such as the ones appearing in \cite{alice-mylene-ptrf, alice-mylene05}) to these implicit master equation representations. Consequently, our technique is simpler, more general and brings into focus the source of the phase transition phenomenon. The underlying methods can and have been adapted to study the extreme singular values and singular vectors of deformations of rectangular random matrices, as well as  the fluctuations \cite{bg-maida-guionnetTCL} and the large deviations \cite{bg-maida-guionnetLDP} of our model.

The paper is organized as follows. In Section \ref{sec:main results}, we state the main results and present the 
integral transforms alluded to above. Section \ref{sec:examples} presents some examples. An outline of the proofs is presented in Section \ref{sec:sketch add}. Exact master equation representations of the eigenvalues and the eigenvectors of the perturbed matrices are derived  in Section \ref{sec:master_equation} and utilized in Section \ref{210709.18h26}   to prove the main results. Technical results needed in these proofs have been relegated to the Appendix.

%% file: mainresults_fin_v6.tex
\section{Main results}\label{sec:main results}

\subsection{Definitions and hypotheses}\label{sec:definitions}
Let  $X_n$ be  an $n \times n$ symmetric (or Hermitian) random matrix whose ordered eigenvalues we denote by $\la_{1}(X_n) \geq \cdots \geq \la_{n}(X_n)$. Let $\mu_{X_{n}}$ be the empirical eigenvalue distribution, \ie, the probability measure defined as
\[
\mu_{X_{n}} = \frac{1}{n} \sum_{j=1}^{n} \delta_{\la_{j}(X_{n})}.
\]
Assume that the probability measure $\mu_{X_{n}}$ converges almost surely weakly, as $n \longrightarrow \infty$, to a non-random compactly supported  \pro measure $\mu_X$. Let $a$ and $b$ be, respectively, the infimum and supremum of the support of $\mu_X$. We suppose  the smallest and largest eigenvalue of $X_n$ converge almost surely to $a$ and $b$.

For a given $r \ge 1$, let $\tta_1\ge\cdots\ge \tta_r$ be deterministic non-zero real numbers, chosen independently of $n$.  For every $n$, let $P_n$ be an $n\times n$ symmetric (or Hermitian) random matrix   having rank $r$ with its $r$ non-zero eigenvalues equal to $\theta_{1}, \ldots, \theta_{r}$.
Let the index $s\in \{0, \ldots, r\}$  be defined such that   $\tta_1\ge\cdots\ge \tta_s>0>\tta_{s+1}\ge\cdots\ge \tta_r$.

Recall that a symmetric (or Hermitian) random matrix is said to be {\it orthogonally invariant} (or {\it unitarily invariant}) if its distribution is invariant under the action of the orthogonal (or unitary) group under conjugation.

We suppose that $X_n$ and $P_n$ are independent and that either $X_n$ or $P_n$ is  orthogonally (or   unitarily) invariant.

\subsection{Notation} Throughout this paper, for $f$ a function and $c\in \R$, we set $$f(c^+):=\lim_{z\downarrow c}f(z)\,;\qquad f(c^-):=\lim_{z\uparrow c}f(z),$$
we also let $\convas$ denote almost sure convergence. The ordered eigenvalues of an $n\times n$ Hermitian matrix $M$ will be denoted by $\la_1(M)\ge\cdots \ge \la_n(M)$. Lastly, for a subspace $F$ of a Euclidian space $E$ and a vector $x\in E$, we denote the norm of the orthogonal projection of $x$ onto $F$ by $\lan x, F\ran$. 

\subsection{Extreme eigenvalues and eigenvectors under additive perturbations}

Consider the rank $r$ additive perturbation of the random matrix $X_n$ given by
$$\wtX = X_n+P_n.$$

\begin{Th}[Eigenvalue phase transition]\label{140709.main}
The extreme eigenvalues of
$\wtX $
exhibit the following    behavior as $n \longrightarrow \infty$. We have that for each $1\le  i \leq s$,
$$\la_i(\wtX) \convas \begin{cases} G_{\mu_X}^{-1}(1/\theta_i)& \textrm{ if } \theta_i>1/G_{\mu_X}(b^+),\\ \\
b &\textrm{ otherwise,}\end{cases}$$
while for each fixed $i>s$, $\la_i(\wtX)\convas b$.

Similarly, for the smallest eigenvalues, we have that for each  $0 \leq j < r-s$,  $$\la_{n-j}(\wtX) \convas
\begin{cases}
G_{\mu_X}^{-1}(1/\theta_{r-j}) & \textrm{ if } \theta_j<1/G_{\mu_X}(a^-),\\ \\
a &    \textrm{ otherwise,}
\end{cases}
$$while for each fixed $j\ge r-s$, $\la_{n-j}(\wtX)\convas a$.

Here,
$$G_{\mu_X}(z)=\int\f{1}{z-t}\ud\mu_X(t) \qquad \textrm{for } z \notin \supp \mu_X,$$
is the Cauchy transform of $\mu_X$,  $G_{\mu_X}^{-1}(\cdot)$ is its functional inverse so that $1/±\pm\infty$ stands for $0$.
\end{Th}

\begin{Th}[Norm of the eigenvector projection]\label{180709.13h39} Consider   ${i_0}\in \{1, \ldots, r\}$  \st $1/\theta_{i_0} \in (G_{\mu_X}(a^-), G_{\mu_X}(b^+))$. For each $n$, define $$\widetilde{\lambda}_{i_0}:=\begin{cases}\la_{i_0}(\wtX)&\textrm{if $\theta_{i_0}>0$,}\\ \\
\la_{n-r+i_0}(\wtX)&\textrm{if $\theta_{i_0}<0$,}\end{cases}$$ and let $\widetilde{u}$ be a unit-norm eigenvector of $\wtX$  associated with the eigenvalue $\widetilde{\lambda}_{i_0}$. Then we have, as $n\lto\infty$,

\flushleft (a) $$|\langle \widetilde{u}, \ker (\theta_{{i_0}}I_n-P_n) \rangle|^{2} \convas  {\f{-1}{\theta_{i_0}^2G_{\mu_X}'(\rho)}}$$
where $\rho=G_{\mu_X}^{-1}(1/\theta_{{i_0}})$ is the limit of $\widetilde{\lambda}_{i_0}$;

\flushleft (b)
$$\langle \widetilde{u}, \oplus_{i\neq {i_0}} \ker (\theta_{i}I_n-P_n) \rangle \convas 0.$$
 \end{Th}

\begin{Th}[Eigenvector phase transition]\label{200709.13h5}When $r=1$, let the sole non-zero eigenvalue of $P_n$  be denoted by $\theta$.
Suppose that
$$ \ff{\theta}\notin (G_{\mu_X}(a^-), G_{\mu_X}(b^+)),\quad\textrm{and}\quad\begin{cases}G_{\mu_X}'(b^+)=-\infty&\textrm{if $\theta>0$,}\\ \\ G_{\mu_X}'(a^-)=-\infty&\textrm{if $\theta<0$.}\end{cases}$$
For each $n$, let $\widetilde{u}$ be a unit-norm eigenvector of $\wtX$ associated with either the largest or smallest eigenvalue depending on whether $\theta>0$ or $\theta<0$, respectively. Then we have $$\lan \widetilde{u}, \ker(\theta I_n-P_n)\ran\convas 0, $$
as $n \longrightarrow \infty$.
\end{Th}

The following proposition allows to assert that in many classical matrix models, such as Wigner or Wishart matrices, the above phase transitions actually occur with a finite threshold. The proposition is phrased in terms of $b$, the supremum of the support of $\mu_X$, but also applies for $a$, the infimum of the support of $\mu_X$. The proof relies on a straightforward computation which we omit.

\begin{propo}[Square-root decay at edge and phase transitions]\label{square root add}
Assume that the limiting eigenvalue distribution $\mu_X$ has a density $f_{\mu_X}$ with a power decay at $b$, i.e., that, as $t\to b$ with $t<b$,  $f_{\mu_X}(t) \sim c (b-t)^\alpha$  for some exponent $\alpha>-1$ and  some constant $c$. Then:
$$ G_{\mu_X}(b^{+}) < \infty\iff \alpha>0  \qquad\textrm{ and } \qquad G'_{\mu_X}(b^{+})=-\infty \iff \al\le 1,$$
so that the phase transitions in Theorems \ref{140709.main} and \ref{200709.13h5} manifest for $\al=1/2$.
\end{propo}

\begin{rmq}[Necessity of eigenvalue repulsion for the eigenvector phase transition]\label{210709.23h19}{\rm 
Under additional hypotheses on the manner in which the empirical eigenvalue distribution of $X_n \convas \mu_X$ as $n \longrightarrow \infty$,  Theorem \ref{180709.13h39} can be generalized to any eigenvalue with limit $\rho$ equal either to  $a$ or $b$ \st $G_{\mu_X}'(\rho)$ is finite. In the same way, Theorem \ref{200709.13h5} can be generalized for any value of $r$. The specific hypothesis has to do with requiring the spacings between the $\la_i(X_n)$'s to be more ``random matrix like'' and exhibit repulsion instead of being ``independent sample like'' with possible clumping. We plan to develop this line of inquiry in a separate paper.}
\end{rmq}

\subsection{Extreme eigenvalues and eigenvectors under multiplicative perturbations}
We maintain the same hypotheses as before so that the limiting probability measure $\mu_X$, the index $s$ and the rank $r$ matrix $P_{n}$ are defined as in Section \ref{sec:definitions}. In addition, we assume that for every $n$, $X_n$ is a non-negative definite matrix and that the limiting probability measure $\mu_X$ 
is not the Dirac mass at zero.

Consider the rank $r$ multiplicative perturbation of the random matrix $X_n$ given by
$$\wtX = X_n\times (I_n+P_n).$$
\begin{Th}[Eigenvalue phase transition]\label{140709.main.multiplicative}
The extreme eigenvalues of
$\wtX$
exhibit the following behavior as $n \longrightarrow \infty$. We have that for $1 \leq i \leq  s$, $$\la_i(\wtX) \convas
\begin{cases}
T_{\mu_X}^{-1}(1/\theta_i)& \textrm{ if }  \theta_i>1/T_{\mu_X}(b^+),\\ \\
b &\textrm{ otherwise,}\end{cases}$$
while for each fixed $i>s$,
$\la_i(\wtX)\convas b$.

In the same way, for the smallest eigenvalues, for each  $0 \leq j < r-s$,
$$\la_{n-r+j}(\wtX) \convas
\begin{cases} T_{\mu_X}^{-1}(1/\theta_j)&\textrm{ if } \theta_j< 1/ T_{\mu_X}(a^-),\\ \\ a &\textrm{ otherwise,}\end{cases}$$
while for each fixed $j\ge r-s$,
$\la_{n-j}(\wtX)\convas a$.

Here,
$$T_{\mu_X}(z)=\int\f{t}{z-t}\ud \mu_X(t)  \qquad \textrm{for } z \notin \supp \mu_X,$$
is the T-transform of $\mu_X$, $T_{\mu_X}^{-1}(\cdot)$ is its functional inverse and $1/\pm \infty$ stands for $0$.
\end{Th}

\begin{Th}[Norm of eigenvector projection]\label{180709.13h39.multiplicative}
Consider  ${{i_0}}\in \{1, \ldots, r\}$  \st $1/\theta_{i_0} \in (T_{\mu_X}(a^-), T_{\mu_X}(b^+))$. For each $n$, define $$\widetilde{\lambda}_{i_0}:=\begin{cases}\la_{i_0}(\wtX)&\textrm{if $\theta_{{i_0}}>0$,}\\ \\
\la_{n-r+{i_0}}(\wtX)&\textrm{if $\theta_{{i_0}}<0$,}\end{cases}$$ and let $\widetilde{u}$ be a unit-norm eigenvector of $\wtX$  associated with the eigenvalue $\widetilde{\lambda}_{i_0}$. Then we have, as $n \longrightarrow \infty$,

\flushleft a)
$$| \langle \widetilde{u}, \ker (\theta_{{i_0}}I_n-P_n) \rangle|^{2} \convas  {\f{-1}{\theta_{{i_0}}^2\rho T_{\mu_X}'(\rho)+\theta_{{i_0}} }},$$
where $\rho=T_{\mu_X}^{-1}(1/\theta_{{i_0}})$ is the limit of $\widetilde{\lambda}_{i_0}$;

\flushleft b)  $$ \langle \widetilde{u}, \oplus_{j\neq {i_0}} \ker (\theta_{j}I_n-P_n) \rangle  \convas 0.$$
\end{Th}

\begin{Th}[Eigenvector phase transition]\label{200709.13h5.multiplicative}When $r=1$, let the sole non-zero eigenvalue of $P_n$  be denoted by $\theta$.  Suppose that
$$ \ff{\theta}\notin (T_{\mu_X}(a^-), T_{\mu_X}(b^+)),\quad\textrm{and}\quad\begin{cases} T_{\mu_X}'(b^+)=-\infty&\textrm{if $\theta>0$,}\\  \\  T_{\mu_X}'(a^-)=-\infty&\textrm{if $\theta<0$.}\end{cases}$$
For each $n$, let $\widetilde{u}$ be the unit-norm eigenvector of $\wtX$ associated with either the largest or smallest eigenvalue depending on whether $\theta>0$ or $\theta<0$, respectively. Then, we have
$$\lan \widetilde{u}, \ker(\theta I_n-P_n)\ran \convas 0$$
as $n \longrightarrow \infty$.
\end{Th}

\begin{propo}[Square-root decay at edge and phase transitions]\label{square root mult}
Assume that the limiting eigenvalue distribution $\mu_X$ has a density $f_{\mu_X}$ with a power decay at $b$ (or $a$ or both), i.e., that, as $t\to b$ with $t<b$,  $f_{\mu_X}(t) \sim c (b-t)^\alpha$  for some exponent $\alpha>-1$ and  some constant $c$. Then:
$$ T_{\mu_X}(b^{+}) < \infty\iff \alpha>0  \qquad\textrm{ and } \qquad T'_{\mu_X}(b^{+})=-\infty \iff \al\le 1,$$
so that the phase transitions in Theorems \ref{140709.main.multiplicative} and \ref{200709.13h5.multiplicative} manifest for $\al=1/2$. \end{propo}

\noindent The analogue of Remark \ref{210709.23h19} also applies here.

\begin{rmq}[Eigenvalues and eigenvectors of a similarity transformation of $X_n$]\label{rem:similar mult}
{\rm
Consider the matrix $S_n = (I_n+P_n)^{1/2}X_n (I_n+P_n)^{1/2}$. The matrices $S_n$ and $\widetilde{X}_n=X_n(I_n+P_n)$ are related by a similarity transformation and hence share the same eigenvalues and consequently the same limiting eigenvalue behavior in Theorem \ref{140709.main.multiplicative}. Additionally, if $\widetilde{u}_i$ is a unit-norm eigenvector of $\widetilde{X}_{n}$ then $\widetilde{w}_i=(I_n+P_n)^{1/2}\widetilde{u}_i$ is an eigenvector of $S_{n}$ and the unit-norm eigenvector $\widetilde{v}_i=\widetilde{w}_i/\|\widetilde{w}_i\|$ satisfies
$$
|\langle \widetilde{v}, \ker(\theta_{i}I_n-P_n) \rangle|^{2} = \frac{ (\theta_{i}+1) |\langle \widetilde{u}, \ker(\theta_{i}I_n-P_n) \rangle|^{2} }{\theta_{i} \langle \widetilde{u}, \ker(\theta_{i} I_n-P_n) \rangle^{2}+1}.
$$
It follows that we obtain the same phase transition behavior and that when $1/\theta_{i} \in (T_{\mu_X}(a^-), T_{\mu_X}(b^+))$,  $$|\langle \widetilde{v}_i, \ker(\theta_{i}I_n-P_n) \rangle|^{2} \convas -\f{\theta_{i}+1}{\theta_{i}T'_{\mu_X}(\rho)}\quad\textrm{ and }\quad \langle \widetilde{v}_i, \oplus_{j\neq i} \ker (\theta_{i}I_n-P_n) \rangle  \convas 0,$$
so that   the analogue of Theorems \ref{180709.13h39.multiplicative} and \ref{200709.13h5.multiplicative} for the eigenvectors of $S_n$ holds.}
\end{rmq}

\subsection{The Cauchy and T transforms in free probability theory}\label{8709.14h}
\subsubsection{The Cauchy transform and its relation to additive free convolution}\label{8709.14h.bigfish}
 The {\it Cauchy transform} of a compactly supported \pro measure $\mu$ on the real line is defined as:  $$G_\mu(z)=\int\f{\ud \mu(t)}{z-t}  \qquad \textrm{for } z \notin \supp \mu.$$
If  $[a,b]$ denotes the convex hull of  the support of $\mu$, then
$$G_\mu(a^-):=\lim_{z\uparrow a}G_\mu(z)\qquad\textrm{  and }\qquad G_\mu(b^+):=\lim_{z\downarrow b}G_\mu(z)$$ exist in $[-\infty, 0)$ and $(0,+\infty]$, respectively and $G_\mu(\cdot)$ realizes decreasing homeomorphisms from $(-\infty, a)$ onto $(G_\mu(a^-),0)$ and from $(b, +\infty)$ onto $(0, G_\mu(b^+))$. Throughout this paper, we shall denote by $G_\mu^{-1}(\cdot)$ the inverses of these homeomorphisms, even though $G_\mu$  can also  define other homeomorphisms on the holes of the support of $\mu$.

The {\it $R$-transform}, defined as $$R_\mu(z):=G_\mu^{-1}(z)-1/z,$$ is the analogue of the logarithm of the Fourier transform for free additive convolution. The free additive convolution of \pro measures on the real line   is denoted by the symbol $\bxp$ and can be characterized as follows.

Let $A_n$ and $B_n$ be independent $n \times n$ symmetric (or Hermitian) random matrices that   are invariant, in law, by conjugation by any orthogonal (or unitary) matrix. Suppose that, as $n \longrightarrow \infty$, $\mu_{A_{n}} \longrightarrow \mu_{A}$ and $\mu_{B_{n}} \longrightarrow \mu_{B}$. Then, free probability theory states that $\mu_{A_n + B_n} \longrightarrow \mu_{A} \bxp \mu_{B}$, a \pro measure which
 can be characterized in terms of the $R$-transform as
$$R_{\mu_A\bxp\mu_B}(z)= R_{\mu_A}(z)+R_{\mu_B}(z).$$ The connection between free additive convolution and $G_\mu^{-1}$ (via the $R$-transform) and the appearance of  $G_\mu^{-1}$ in Theorem \ref{140709.main} could be of independent interest to free probabilists.

\subsubsection{The $T$-transform and its relation to multiplicative free convolution}\label{subsection: T-transform}
In the case where $\mu\neq \delta_0$ and the support of $\mu$ is contained in $[0, +\infty)$,
 one   also defines its {\it $T$-transform}
 $$T_\mu(z)=\int\f{t}{z-t}\ud \mu(t)  \qquad \textrm{for } z \notin \supp \mu_X,$$
which realizes decreasing homeomorphisms from $(-\infty, a)$ onto $(T_\mu(a^-),0)$ and from $(b, +\infty)$ onto $(0, T_\mu(b^+))$. Throughout this paper, we shall denote by $T_\mu^{-1}$ the inverses of these homeomorphisms, even though $T_\mu$   can also  define other  homeomorphisms on the holes of the support of $\mu$.

The {\it $S$-transform}, defined as
$$S_\mu(z):=(1+z)/({z} {T_\mu^{-1}(z)}),$$
is the analogue  of the Fourier transform for  free multiplicative convolution $\bxt$. The free multiplicative convolution of two \pro measures $\mu_A$ and $\mu_B$ is denoted by the symbols $\bxt$ and can be characterized as follows.

Let $A_n$ and $B_n$ be independent $n \times n$ symmetric (or Hermitian) positive-definite random matrices that   are invariant, in law, by conjugation by any orthogonal (or unitary) matrix. Suppose that, as $n \longrightarrow \infty$, $\mu_{A_{n}} \longrightarrow \mu_{A}$ and $\mu_{B_{n}} \longrightarrow \mu_{B}$. Then, free probability theory states that $\mu_{A_n \cdot B_n} \longrightarrow \mu_{A} \bxt \mu_{B}$, a \pro measure which can be characterized in terms of the $S$-transform as
$$S_{\mu_A\bxt\mu_B}(z)= S_{\mu_A}(z)S_{\mu_B}(z).$$ The connection between free multiplicative convolution and $T_\mu^{-1}$ (via the $S$-transform)  and the 
appearance of  $T_\mu^{-1}$ in Theorem \ref{140709.main.multiplicative} could be of independent interest to free probabilists.


\subsection{Extensions}

\begin{rmq}[Phase transition in non-extreme eigenvalues]\label{220709.12h27}{\rm Theorem \ref{140709.main} can easily be adapted to describe the phase transition in the eigenvalues of  $X_n+P_n$ which fall in the ``holes'' of the support of $\mu_X$. Consider $c<d$ \st almost surely, for $n$ large enough, $X_n$ has no eigenvalue in the interval $(c,d)$. It  implies  that $G_{\mu_X}$ induces a decreasing  homeomorphism, that we shall denote by $G_{\mu_X,(c,d)}$, from  the interval $(c,d)$  onto  the interval $(G_{\mu_X}(d^-), G_{\mu_X}(c^+))$. Then it can be proved that almost surely, for $n$ large enough, $X_n+P_n$ has no eigenvalue in the interval $(c,d)$, except  if some of the $1/\theta_i$'s are in the interval  $(G_{\mu_X}(d^-),G_{\mu_X}(c^+))$, in which case for each such index $i$, one  eigenvalue of $X_n+P_n$ has limit $ G_{\mu_X,(c,d)}^{-1}(1/\theta_i)$ as $n \longrightarrow \infty$.}
\end{rmq}

\begin{rmq}[Isolated eigenvalues  of $X_n$  outside the support of $\mu_X$]\label{210709.18h27}{\rm Theorem \ref{140709.main} can also easily be adapted to the case where
 $X_n$ itself has  isolated eigenvalues in the sense that some of its eigenvalues have limits out of the support of $\mu_X$. More formally, let us replace the assumption that the smallest and largest eigenvalues of $X_n$ tend to the infimum $a$ and the supremum $b$ of the support of $\mu_X$ by the following one.
 \begin{quote}
 \emph{There exists some real numbers
 $$\ell_1^+, \ldots, \ell_{p+}^+\in (b, +\infty) \qquad \textrm{and }\qquad  \ell_{1}^-,\ldots, \ell_{p^-}^-\in (-\infty,a)
$$
\st for all $1\le j\le p^+$, $$\la_j(X_n)\convas \ell_j^+$$
and for all $1\le j\le p^-$, $$\la_{n+1-j}(X_n)\convas \ell_j^-.$$ Moreover, $\la_{1+p^+}(X_n)\convas b$ and $\la_{n-(1+p^-)}(X_n)\convas a$.}\end{quote}
Then an attentive look at the proof of Theorem \ref{140709.main} shows that it still holds, in the following sense (we only present the point of view of the largest eigenvalues): the matrix  $\wtX$ still has eigenvalues tending to the $\ell_j^+$'s, but also,
for each $1\le i\le s$ \st $\tta_i>1/G_{\mu_X}(b^+)$, one  eigenvalue tending to  $G_{\mu_X}^{-1}(1/\tta_i)$, all other largest eigenvalues tending to $b$.}\end{rmq}

   \begin{rmq}[Other perturbations of  $X_n$]{\rm The previous remark forms the basis for an iterative application of our theorems to other perturbational models, such as $\widetilde{X}=\sqrt{X}(I+P)\sqrt{X}+Q$ for example. Another way to deal with such perturbations is to first derive the corresponding {\it master equations} representations that describe how the eigenvalues and eigenvectors of $\widetilde{X}$  are related to the eigenvalues and eigenvectors of $X$ and the perturbing matrices, along the lines of  Proposition \ref{9709.10h} for additive or multiplicative perturbations of Hermitian matrices.}\end{rmq}

\begin{rmq}[Random matrices with Haar-like eigenvectors]\label{rem:haarlike extension}
{\rm
Let $G$ be an $n \times m$ Gaussian random matrix with independent real (or complex) entries that are normally distributed with mean $0$ and variance $1$. Then the matrix $X = GG^{*}/m$ is orthogonally (or unitarily) invariant. Hence one can choose an orthonormal basis $(U_1,\ldots, U_n)$ of eigenvectors of $X$ \st the matrix $U$ with columns $U_1, \ldots, U_n$ is Haar-distributed.    When $G$ is a Gaussian-like matrix, in the sense that its entries are i.i.d. with mean zero and variance one, then upon placing adequate restrictions on the higher order moments,  for non-random unit norm vector $x_n$, the vector $U^{*}x_{n}$ will be close to uniformly distributed on the unit real (or complex) sphere \cite{s84,s89,s90}. Since our proofs rely heavily on the properties of unit norm vectors uniformly distributed on the $n$-sphere,  they could possibly be adapted to the setting where the unit norm vectors are close to uniformly distributed.}
\end{rmq}

\begin{rmq}[Setting where eigenvalues of $P_n$ are not fixed]\label{rem:random Pn}
{\rm
Suppose that $P_n$ is a random matrix independent of $X_n$, with exactly $r$ non-zero eigenvalues   given by $\theta_1^{(n)},\ldots, \theta_r^{(n)}$. Let $\theta_{i}^{(n)} \convas \theta_{i}$ as $n \longrightarrow \infty$. Using \cite[Cor. 6.3.8]{hj85} as in Section \ref{71009.15h19}, one can easily see that our results will also apply in this case.}\end{rmq}

\noindent The analogues of Remarks   \ref{220709.12h27}, \ref{210709.18h27}, \ref{rem:haarlike extension} and \ref{rem:random Pn} for the multiplicative setting also hold here. 
In particular, Wishart matrices with $c>1$ (cf Section \ref{subsectionWishart}) gives an illustration of the case where there is a hole in the support of $\mu_X$.

%% file: examples_fin_v6.tex
\section{Examples}\label{sec:examples}
We now illustrate our results with some concrete computations. The key to applying our results lies in being able to compute the Cauchy or $T$ transforms of the probability measure $\mu_X$ and their associated functional inverses. In what follows, we focus on settings where the transforms and their inverses can be expressed in closed form. In settings where the transforms are algebraic so that they can be represented as solutions of polynomial equations, the techniques and software developed in \cite{re08} can be utilized. In more complicated settings, one will have to resort to numerical techniques.

\subsection{Additive perturbation of a Gaussian Wigner matrix}\label{sec:gaussian square}
Let $X_n$ be an $n\times n$ symmetric (or Hermitian) matrix with independent, zero mean, normally distributed entries with variance $\sigma^2/n$ on the diagonal and $\sigma^2/(2n)$ on the off diagonal. It is known that the spectral measure of $X_n$ converges almost surely to the famous semi-circle distribution with density
$$
\ud \mu_{X}(x)= \frac{\sqrt{4\sigma^{2} - x^{2}}}{2\sigma^{2}\pi} \ud x \qquad \textrm{ for } x \in [-2\sigma, 2\sigma].
$$
It is known that the extreme eigenvalues converge almost surely to the endpoints of the support \cite{agz09}. Associated with the spectral measure, we have  $$G_{\mu_{X}}(z)=\f{z-\operatorname{sgn}(z)\sqrt{z^2-4\sigma^2}}{2\sigma^2},\quad \textrm{ for } z\in (-\infty,-2\sigma)\cup(2\sigma,+\infty), $$ $G_{\mu_X}(\pm 2\sigma)=\pm\sigma$ and $G_{\mu_X}^{-1}(1/\theta)=\theta+\f{\sigma^2}{\theta}$.

Thus for a $P_n$ with $r$ non-zero eigenvalues   $\theta_1\geq \cdots\geq\theta_s>0>\theta_{s+1}\geq  \cdots\geq\theta_r$, by Theorem \ref{140709.main}, we have for $1 \leq i \leq s$,
\begin{equation*}\label{eq:ex1 eigval}
\la_i(X_n+P_n)  \convas
\begin{cases}
\theta_i+\f{\sigma^2}{\theta_i} & \textrm{ if } \theta_i > \sigma \\
2\sigma   & \textrm{ otherwise},
\end{cases}
\end{equation*}
as $n \longrightarrow \infty$.  This result has already been established in \cite{pf07} for the symmetric case and in \cite{sp06} for the Hermitian case. Remark \ref{rem:haarlike extension} explains why our results should hold for Wigner matrices of the sort considered in \cite{sp06,pf07}. 

In the setting where $r=1$ and  $P=\theta \,uu^{*}$, let $\widetilde{u}$  be a unit-norm eigenvector of $X_n+P_n$  associated with its largest eigenvalue. By Theorems \ref{180709.13h39} and \ref{200709.13h5}, we have
\begin{equation*}\label{eq:ex1 eigvec}
|\lan \widetilde{u}, u \ran|^2\convas \begin{cases}
1-\f{\sigma^2}{\theta^2}&\textrm{if $\theta\geq \sigma$,}\\
0&\textrm{if $\theta< \sigma$.}
\end{cases}
\end{equation*}

\subsection{Multiplicative perturbation of a random Wishart matrix}\label{subsectionWishart}
Let $G_n$ be an $n\times m$ real (or complex) matrix with independent, zero mean, normally distributed entries with variance $1$. Let $X_n = G_nG_n^{*}/m$. It is known \cite{mp67} that, as $n,m \longrightarrow \infty$ with $n/m \to c>0$,  the spectral measure of $X_n$  converges almost surely to the famous Mar\v{c}enko-Pastur distribution with density
$$\ud \mu_{X}(x):=\f{1}{2\pi cx}\sqrt{(b-x)(x-a)}\one_{[a,b]}(x)\ud x + \max\left(0,1-\frac{1}{c}\right) \delta_{0},$$
where $a=  (1-\sqrt{c})^2$ and $b= (1+\sqrt{c})^2$. It is known \cite{book-bai-silverstein} that the extreme eigenvalues converge almost surely to the endpoints of this support.

Associated with this spectral measure, we have
\beq
T_{\mu_X}(z)&=&\f{z-c-1-\operatorname{sgn}(z-a) \sqrt{(z-a)(z-b)}}{2c} 
\eeq
$T_{\mu_X}(b^+)=1/\sqrt{c}$, $T_{\mu_X}(a^-)\;=\;-1/\sqrt{c}$ and
\beq
T_{\mu_X}^{-1}(z)=\f{(z+1)(cz+1)}{z}.
\eeq

When $c>1$, there is an atom at zero so that the smallest eigenvalue of $X_n$ is identically zero. For simplicity, let us   consider the setting when $c <1$ so that the extreme eigenvalues of $X_n$ converge almost surely to $a$ and $b$. Thus for $P_n$ with $r$ non-zero eigenvalues  $\theta_1\geq \cdots\geq\theta_s>0>\theta_{s+1}\geq  \cdots\geq\theta_r$, with $l_{i} :=\theta_{i} + 1$,
 for $c<1$, by Theorem  \ref{140709.main.multiplicative}, we have for $1\le i\le s$,
\[
\la_i(X_n(I_n+P_n))  \convas
\begin{cases}
l_i\lf(1+\f{c}{l_i-1}\ri) & \textrm{ if  } |l_i -1| > \sqrt{c} \\
b  & \textrm{otherwise},
\end{cases}
\]
as $n \longrightarrow \infty$. An analogous result for the smallest eigenvalue may be similarly derived by making the appropriate substitution for $a$ in Theorem  \ref{140709.main.multiplicative}.
Consider the matrix $S_n = (I_n+P_n)^{1/2}X_n (I_n+P_n)^{1/2}$. The matrix $S_n$ may be interpreted as a Wishart distributed sample covariance matrix with ``spiked'' covariance $I_n + P_n$. By Remark \ref{rem:similar mult}, the above result applies for the eigenvalues of $S_n$ as well.  This result for the largest eigenvalue of spiked sample covariance matrices was established in \cite{bbp05,p07} and for the extreme eigenvalues in \cite{bs06}.

In the setting where $r=1$ and  $P=\theta \,uu^{*}$, let $l= \theta+1$ and let $\widetilde{u}$   be a unit-norm eigenvector of $X_n(I+P_n)$  associated with its largest (or  smallest, depending on whether $l>1$ or $l<1$)  eigenvalue. By Theorem \ref{200709.13h5.multiplicative}, we have
\bes\label{240909.14h22}|\lan \widetilde{u}, u \ran|^2\convas \begin{cases}
\f{ (l-1)^2-c}{(l-1)[c(l+1)+l-1]}&\textrm{if $|l-1|\geq \sqrt{c}$,}
\\
0&\textrm{if $|l-1|< \sqrt{c}$.}
\end{cases}\ees

Let $\widetilde{v}$ be a unit eigenvector of $S_n = (I_n+P_n)^{1/2}X_n (I_n+P_n)^{1/2}$ associated with its largest (or  smallest, depending on whether $l>1$ or $l<1$)  eigenvalue. Then, by Theorem \ref{200709.13h5.multiplicative} and Remark \ref{rem:similar mult}, we have
\bes\label{240909.14h222}|\lan \widetilde{v}, u \ran|^2\convas \begin{cases}
\f{ 1-\f{c}{(l-1)^2}}{1+\f{c}{l-1}}&\textrm{if $|l-1|\geq \sqrt{c}$,}
\\
0&\textrm{if $|l-1|< \sqrt{c}$.}
\end{cases}\ees
The result has been established in \cite{p07} for the eigenvector associated with the  largest eigenvalue. We generalize it to the eigenvector associated with the smallest one. 

We note that symmetry considerations imply that when $X$ is a Wigner matrix then $-X$ is a Wigner matrix as well. Thus an analytical characterization of the largest eigenvalue of a Wigner matrix directly yields a characterization of the smallest eigenvalue as well. This trick cannot be applied for Wishart matrices since Wishart matrices do not exhibit the symmetries of Wigner matrices. Consequently, the smallest and largest eigenvalues and their associated eigenvectors of Wishart matrices have to be treated separately. Our results facilitate such a characterization.

%% file: sketch-proofs_fin_v6.tex
\section{Outline of the proofs}\label{sec:sketch add}
We now provide an outline of the proofs. We focus on Theorems \ref{140709.main}, \ref{180709.13h39} and \ref{200709.13h5},  which describe   the phase transition in the extreme eigenvalues and associated eigenvectors of $X+P$ (the index $n$ in $X_n$ and $P_n$ has been suppressed for brevity). An analogous argument applies for the multiplicative perturbation setting.

Consider the setting where $r=1$, so that $P = \theta \,uu^{*}$, with $u$ being a unit norm column vector.  Since either $X$ or $P$ is assumed to be invariant, in law, under  orthogonal (or unitary) conjugation, one can, without loss of generality, suppose that $X=\diag(\la_1, \ldots, \la_n)$ and that  $u$ is uniformly distributed on the unit $n$-sphere.

\subsection{Largest eigenvalue phase transition}\label{sec:max eig sketch}
The eigenvalues of $X+P$ are the solutions of the equation
$$ \det(zI - (X+P)) = 0 .$$
Equivalently, for $z $ so that $zI-X$ is invertible, we have $$zI-(X+P)=(zI-X) \cdot (I-(zI-X)^{-1}P),$$ so that
$$ \det(zI - (X+P)) = \det(zI - X) \cdot \det (I - (zI -X)^{-1}P).$$
Consequently, a simple argument reveals that the $z$ is an eigenvalue of $X+P$ and not an eigenvalue of $X$ if and only if $1$ is an eigenvalue of the matrix $(zI-X)^{-1}P$. But $(zI-X)^{-1}P= (zI-X)^{-1}\theta \,uu^{*}$ has rank one, so its only non-zero eigenvalue will equal its trace, which in turn is equal to $\theta G_{\mu_n}(z)$, where ${\mu_n}$ is a ``weighted" spectral measure of $X$, defined by \begin{equation*}\label{eq:master eq mur1}\qquad
\mu_n = \sum_{k=1}^{n} | u_{k} |^{2} \delta_{\la_k} \qquad \textrm{ (the $u_k$'s are the coordinates of $u$).}
\end{equation*}
Thus any $z$ outside the  spectrum of $X$  is an eigenvalue of $X+P$ if and only if
\begin{equation}\label{eq:master eq r1}
 \sum_{k=1}^{n}\dfrac{| u_{k} |^{2}}{z-\lambda_{k}}=:G_{\mu_n}(z)  = \ff{\theta}.
\end{equation}
Equation (\ref{eq:master eq r1}) describes the relationship between the eigenvalues of $X+P$ and the eigenvalues of $X$ and the dependence on the coordinates of the vector $u$ (via the measure $\mu_n$).

This is where randomization simplifies analysis. Since $u$ is a \text{random} vector with uniform distribution on the unit $n$-sphere, we have that for large $n$, $|u_{k}|^{2} \approx \frac{1}{n}$ with high probability. Consequently, we have  $\mu_n\approx \mu_{X}$ so that $G_{\mu_n}(z) \approx G_{\mu_{X}}(z)$. Inverting equation (\ref{eq:master eq r1}) after substituting these approximations yields the location of the largest  eigenvalue to be $G_{\mu_{X}}^{-1}(1/\theta)$ as in Theorem \ref{140709.main}.

The phase transition for the extreme eigenvalues emerges because under our assumption that  the limiting probability measure $\mu_X$ is compactly supported on $[a,b]$, the Cauchy transform $G_{\mu_{X}}$ is defined \textit{outside} $[a,b]$ and unlike what happens for $G_{\mu_n}$, we do not always have $G_{\mu_X}(b^+)=+\infty$. Consequently, when $1/\theta <G_{\mu_X}(b^+)$, we have that $\lambda_{1}(\widetilde{X}) \approx G_{\mu_X}^{-1}(1/\theta)$ as before.  However, when $1/\theta \geq G_{\mu_X}(b^+)$, the phase transition manifests and $\lambda_{1}(\widetilde{X}) \approx \lambda_{1}(X) = b$.

An extension of these arguments for fixed $r>1$ yields the general result and constitutes the most transparent justification, as sought by the authors in \cite{bbp05}, for the emergence of this phase transition phenomenon in such perturbed random matrix models.  We rely on concentration inequalities to make the arguments rigorous.

\subsection{Eigenvectors phase transition}\label{250909.17h06}Let $\widetilde{u}$ be a unit eigenvector of $X+P$ associated with the eigenvalue $ z $ that satisfies (\ref{eq:master eq r1}).  From the relationship $(X+P)\widetilde{u}=z \widetilde{u}$, we deduce that, for $P=\theta \,uu^{*}$,
$$\qquad\qquad
(zI-X)\widetilde{u}=P\widetilde{u}
=\theta uu^{*}\widetilde{u}
= (\theta u^{*}\widetilde{u}).u\qquad\qquad\textrm{(because $u^{*}\widetilde{u}$ is a scalar)},$$ implying that $\widetilde{u}$ is proportional to $(zI-X)^{-1}u$.

Since $\widetilde{u}$ has unit-norm, \be\label{250909.17h12}\widetilde{u}=\f{(zI-X)^{-1}u}{\sqrt{u^{*}(zI-X)^{-2}u}}\ee and \be\label{250909.17h05}|\lan \widetilde{u}, \ker(\theta I-P)\ran|^2=|u^{*}\widetilde{u}|^2=\f{(u^{*}(zI-X)^{-1}u)^2}{u^{*}(zI-X)^{-2}u}=\f{G_{\mu_n}(z)^2}{\int  \f{\ud \mu_n(t)}{(z-t)^2}}=\ff{\theta^2\int  \f{\ud \mu_n(t)}{(z-t)^2}}.\ee
Equation \eqref{250909.17h12} describes the relationship between the eigenvectors of $X+P$ and the eigenvalues of $X$ and the dependence on the coordinates of the vector $u$ (via the measure $\mu_n$).

Here too, randomization simplifies analysis since for large $n$, we have $\mu_n\approx \mu_{X}$ and $z\approx \rho$. Consequently, $$\int  \f{\ud \mu_n(t)}{(z-t)^2}\approx \int  \f{\ud \mu_X(t)}{(\rho-t)^2}=-G_{\mu_X}'(\rho),$$
so that
when $1/\theta <G_{\mu_X}(b^+)$, which implies that $\rho>b$,    we have
 $$|\lan \widetilde{u}, \ker(\theta I-P)\ran|^2\convas
 \ff{\theta^2\int  \f{\ud \mu_X(t)}{(\rho-t)^2}}=\f{-1}{\theta^2G_{\mu_X}'(\rho)}>0,$$
 whereas when $1/\theta \ge G_{\mu_X}(b^+)$ and $G_{\mu_X}$ has infinite derivative at $\rho=b$, we have $$\lan \widetilde{u}, \ker(\theta I-P)\ran\convas
0. $$
An extension of these arguments for fixed $r>1$ yields the general result and brings into  focus the connection between the eigenvalue phase transition and the associated eigenvector phase transition. As before, concentration inequalities allow us to make these arguments rigorous.

%% file: masterequations_fin_v6.tex
\section{The exact master equations for the perturbed eigenvalues and eigenvectors}\label{sec:master_equation}
In this section, we provide the $r$-dimensional analogues of the master equations \eqref{eq:master eq r1} and \eqref{250909.17h12} employed in our outline of the proof.

\begin{propo}\label{9709.10h}Let us fix some positive integers $1\le r\le n$. Let $X_n=\diag(\la_1, \ldots, \la_n)$ be a diagonal $n\times n$  matrix and $P_n=U_{n,r}\Theta U_{n,r}^{*}$, with $\Theta=\diag(\tta_1,†\ldots, \tta_r)$ an $r\times r$ diagonal matrix  and $U_{n,r}$ an $n\times r$ matrix with orthonormal columns, \ie,  $U_{n,r}^{*}U_{n,r}=I_r$).

a) Then any $z\notin\{\la_1,†\ldots, \la_n\}$ is an eigenvalue of ${\widetilde{X}}_n:=X_n+P_n$  if and only if  the $r\times r$ matrix \be\label{130410.1flo}M_{n}(z):=I_r-U_{n,r}^{*}(zI_n-X_n)^{-1}U_{n,r}\Theta\ee is singular. In this case, \be\label{8.10.10.19h55}\dim \ker(zI_n-{\widetilde{X}}_n) =\dim\ker M_n(z)\ee   and
for all $x\in \ker(zI_n-{\widetilde{X}}_{n})$, we have $U_{n,r}^{*}x \in \ker M_{n}(z)$   and \be\label{130410.1flo-zurich}x=(zI_n-X_{n})^{-1}U_{n,r}\Theta U_{n,r}^{*}x.\ee


b) Let $u^{(n)}_{k,l}$ denote the $(k,l)-$th element of the $n \times r$ matrix $U_{n,r}$ for $k = 1, \ldots, n$ and $l = 1, \ldots, r$. Then for all $i,j=1,\ldots, r$, the $(i, j)$-th entry of  the  matrix $I_r-U_{n,r}^{*}(zI_n-X_n)^{-1}U_{n,r}\Theta$ can be expressed as \begin{equation}\label{eq:master G}\one_{i=j}-\theta_jG_{\mu_{i,j}^{(n)}}(z),\end{equation} where $\mu_{i,j}^{(n)}$ is the complex measure defined by  $$\mu_{i,j}^{(n)}=\sum_{k=1}^n \ovl{u}_{k,i}^{(n)}u_{k,j}^{(n)}\delta_{\la_k}$$ and $G_{\mu_{i,j}^{(n)}}$ is the Cauchy  transform of $\mu_{i,j}^{(n)}$.

c) In the setting where $\widetilde{X} = X_n\times (I_n+P_n) $ and $P_n = U_{n,r} \Theta U_{n,r}^{*}$ as before, we obtain the analog of a) by replacing every occurrence, in \eqref{130410.1flo} and \eqref{130410.1flo-zurich}, of $(zI_n-X_n)^{-1}$ with $(zI_n-X_n)^{-1}X_n$. We obtain the analog of b) by replacing the Cauchy transform in (\ref{eq:master G}) with the $T$-transform.
\end{propo}

\begin{pr}Part a) is proved, for example, in \cite[Th. 2.3]{agg88}. 
Part b) follows  from a straightforward computation of the $(i,j)$-th entry of  $U_{n,r}^*(zI_n-X_{n})^{-1}U_{n,r}\Theta$.
Part c) can be proved in the same way.\end{pr}




%% file: proofsquareadd_fin_v6.tex
\section{Proof of Theorem \ref{140709.main}}\label{210709.18h26}


The sequence of steps described below yields the desired proof:
\begin{enumerate}
\item{The first, rather trivial, step in the proof of Theorem \ref{140709.main} is to use Weyl's interlacing inequalities to prove that any  extreme eigenvalue of $\wtX$ which does not tend to a limit in $\R\bck[a,b]$  tends either to $a$ or $b$.}
\item{Then, we utilize the  ``master equations" of Section \ref{sec:master_equation}   to express the extreme eigenvalues of $\wtX$ as the $z$'s \st   a certain random $r\times r$  matrix $M_n(z)$ is singular.}
\item{We then exploit convergence properties of certain analytical functions (derived in the appendix) to prove that almost surely,  $M_n(z)$ converges  to a certain diagonal matrix $M_{G_{\mu_X}}(z)$, uniformly in $z$.}
\item{We then invoke a continuity lemma (see Lemma \ref{140709.19h26} - derived next) to claim that almost surely, the $z$'s \st $M_n(z)$ is singular ({\it i.e.} the extreme eigenvalues of $\wtX$) converge to the $z$'s \st $M_{G_{\mu_X}}(z)$ is singular.}
\item{We conclude  the proof by noting that, for our setting, the $z$'s \st $M_{G_{\mu_X}}(z)$ is singular are precisely the $z$'s \st for some $i\in \{1, \ldots, r\}$, $G_{\mu_X}(z)=\ff{\tta_i}$.  Part (ii) of Lemma \ref{140709.19h26}, about the rank of $M_n(z)$, will be useful to assert that when the $\tta_i$'s are pairwise distinct, the multiplicities of the isolated eigenvalues  are all equal to one.}
\end{enumerate}

\subsection{A  continuity lemma for the zeros of certain analytic functions}

We now prove a continuity lemma that will be used in the proof of Theorem \ref{140709.main}. We note that nothing in its hypotheses is random. As hinted earlier, we will invoke it to localize the extreme eigenvalues of $\wtX$.

For $z\in \C$ and $E$ a closed subset of $\R$, set $\operatorname{d}(z, E)=\min_{x\in E}|z-x|$.
\begin{lem}\label{140709.19h26}Let us fix a positive integer $r$, a family $\theta_1,†\ldots, \theta_r$ of pairwise distinct nonzero real numbers, two real numbers  $a<b$, an analytic  function $G(z)$ of the variable $z \in \C\bck[a,b]$   \st
\begin{itemize}
\item[a)] $G(z)\in \R  \iff z \in \R$,
\item[]
\item[b)]  for all  $z\in \R\bck[a,b]$,   $G'(z)<0$,
\item[]
\item[c)] $G(z) \longrightarrow 0$ as $|z| \longrightarrow\infty$.
\end{itemize}
Let us define, for $z\in \C\bck[a,b]$, the $r\times r$ matrix \be\label{220709.22h52}M_G(z)=\Diag( 1-\theta_1G(z),\ldots, 1-\theta_rG(z)),\ee and denote by $z_1>\cdots > z_p$  the $z$'s \st $M_G(z)$ is singular, where $p\in\{0,\ldots, r\}$ is identically equal to the number of $i$'s \st $G(a^-)<1/\tta_i<G(b^+)$.

Let us also consider two sequences $a_n$, $b_n$ with respective limits $a$, $b$ and, for each $n$, a function $M_{n}(z)$, defined on $z \in \C\bck[a_n, b_n]$, with values  in the set of $r\times r$ complex matrices \st the  entries of $M_n(z) $ are analytic functions of $z$. We suppose that
\begin{itemize}
\item[d)] for all $n$, for all $z\in \C\bck\R$, the matrix $M_{n}(z)$ is invertible,
\item[]
\item[e)] $M_{n}(z)$ converges, as $n \longrightarrow \infty$, to the function $M_G(z)$, uniformly on $\{z\in \C\ste \operatorname{d}(z, [a,b])\geq \eta \}$, for all $\eta >0$.
\end{itemize}
Then

\begin{itemize}
\item[(i)] there exists $p$ real sequences $z_{n,1}> \ldots > z_{n,p}$ converging respectively to $z_1, \ldots, z_p$   \st for any  $\eps>0$ small enough, for   $n$ large enough,   the $z$'s in $\R\bck[a-\eps, b+\eps]$ \st $M_n(z)$ is singular are exactly $z_{n,1}, \ldots, z_{n,p}$,
\item[]
\item[(ii)] for $n$ large enough, for each $i$, $M_n(z_{n,i})$ has rank $r-1$.\end{itemize}
\end{lem}

\begin{pr}
Note firstly that the $z$'s \st $M_G(z)$ is singular are the $z$'s \st for a certain $j\in \{1,†\ldots, r\}$, \be\label{8.10.2010.1}1-\tta_jG(z)=0.\ee Since the $\tta_j$'s are pairwise distinct, for any $z$, there cannot exist more than one $j\in \{1,€ \ldots, r\}$ \st \eqref{8.10.2010.1} holds. As a consequence, for all $z$, the rank of $M_G(z)$ is either $r$ or $r-1$. Since the set of matrices with rank at least $r-1$ is open in the set of $r\times r$ matrices, once (i) will be proved, (ii) will follow.

Let us now prove (i). Note firstly that by c), there exists $R>\max\{|a|,|b|\}$ \st for $z $ \st $|z|\geq R$, $|G(z)|\leq \min_i \ff{2|\theta_i|}$.
For any such $z$, $|\det M_G(z)|>2^{-r}$.
By  e),  it follows  that for $n$ large enough,   the $z$'s \st $M_{n}(z)$ is singular satisfy $|z|>R$. By d), it even follows that the   $z$'s \st $M_{n}(z)$ is singular   satisfy $z\in[-R,R]$.

Now, to prove (i),  it suffices to prove that for all $c,d\in \R\bck([a,b]\cup \{z_1, \ldots, z_p\})$ \st $c<d<a$ or $b<c<d$, we have that as $n \longrightarrow \infty$:
\begin{itemize}\item[(H)\quad]
the number of  $z$'s in $(c,d)$ \st $\det M_{n}(z)=0$,  denoted by   $\operatorname{Card}_{c,d}(n)$ tends to $\operatorname{Card}_{c,d}$,  the cardinality of the $i$'s in $\{1, \ldots, p\}$ \st $c<z_i<d$.
\end{itemize}

To prove (H), by additivity, one can suppose that $c$ and $d$ are close enough to have $\operatorname{Card}_{c,d}= 0$ or $1$.
Let us define $\gamma$ to be the circle with diameter $[c,d]$. By a) and since $c,d\notin \{z_1, \ldots, z_p\}$, $\det M_G(\cdot)$ does not vanish on $\gamma$, thus  $$\operatorname{Card}_{c,d}=\ff{2i\pi}\int_\gamma \f{\partial_z \det M_G(z)}{\det M_G(z)}\ud z= \lim_{n\to\infty}\ff{2i\pi}\int_\gamma \f{\partial_z \det M_{n}(z)}{\det M_{n}(z)}\ud z,$$ the last equality following from e).
   It follows that for $n$ large enough, $\operatorname{Card}_{c,d}(n)=\operatorname{Card}_{c,d}$ (note that since  $\operatorname{Card}_{c,d}= 0$ or $1$, no ambiguity due to the orders of the zeros has to be taken into account here).
   \end{pr}

\subsection{Proof of Theorem \ref{140709.main}}

\subsubsection{First step: consequences of Weyl's interlacing inequalities}
Note  that  Weyl's interlacing inequalities imply that for all $1\le i\le n$,
\be\label{weyl-261010}\la_{i+(r-s)}(X_n)\le \la_i(\wtX)\le \la_{i-s}(X_n),\ee where we employ the convention that $\la_k(X_n)=-\infty$ is $k>n$ and $+\infty$ if $k\le 0$.
It follows that the empirical spectral measure of $\wtX \convas \mu_X$ because the empirical spectral measure of $X_n$ does as well.

Since $a$ and $b$ belong to the support of $\mu_X$, we have, for all  $i\ge 1$ fixed, $$\liminf_{n\to\infty} \la_i(X_n)\ge b\qquad\textrm{ and }\qquad \limsup_{n\to\infty} \limsup_{n\to\infty} \la_{n+1-i}(X_n)\le a.$$ By the hypotheses that $$\la_1(X_n)\convas b\qquad\textrm{ and }\qquad  \la_n(X_n)\convas a,$$  it follows that for all  $i\ge 1$ fixed, $\la_i(X_n)\convas b$ and $\la_{n+1-i}(X_n)\convas a$.

By \eqref{weyl-261010}, we deduce both following relation \eqref{1342010.bene} and \eqref{27102010.10h}: for all $i\ge 1$ fixed, we have  \be\label{1342010.bene}\liminf_{n\to\infty} \la_i(\wtX)\ge b\qquad\textrm{ and }\qquad \limsup_{n\to\infty} \la_{n+1-i}(\wtX)\le a\ee and for all $i>s$ (resp. $i\ge r-s$) fixed, we have \be\label{27102010.10h}\la_i(\wtX)\convas b\qquad \textrm{ (resp. } \la_{n-i}(\wtX)\convas a\textrm{).}\ee
Equations \eqref{1342010.bene} and \eqref{27102010.10h} are the first step in the proof of Theorem \ref{140709.main}. They   state that any  extreme eigenvalue of $\wtX$ which does not tend to a limit in $\R\bck[a,b]$  tends either to $a$ or $b$. Let us now prove the crux of the theorem related to the isolated eigenvalues.

\subsubsection{Isolated eigenvalues in the setting where $\theta_1, \ldots, \tta_r$ are pairwise distinct}\label{section18.09.09.10h35}
   In this section, we assume that the eigenvalues $\tta_1, \ldots, \tta_r$ of the perturbing matrix $P_n$ to be pairwise distinct. In the next section, we shall remove this hypothesis by an approximation process.

For a momentarily fixed  $n$,  let the eigenvalues of $X_n$ be denoted by $\la_1\ge  \ldots \ge \la_n$. Consider orthogonal (or unitary) $n\times n$ matrices $U_X$, $U_P$ that diagonalize $X_n$ and $P_n$, respectively, \st  \bes\label{140709.19h47}X_n=U_X \Diag(\la_1,\ldots, \la_n)U_X^{*}\, ,\qquad  P_n=U_P\Diag(\theta_1,\ldots, \theta_r, 0, \ldots, 0)U_P^{*}.\ees

The spectrum of $X_n+P_n$ is  identical to the  spectrum  of the matrix
\begin{equation}\label{eq:the equivalent matrix 1}
\Diag(\la_1,\ldots, \la_n)+\underbrace{U_X^{*}U_P}_{\textrm{denoted by }U_n}\Diag(\theta_1,\ldots, \theta_r, 0, \ldots, 0)U_P^{*}U_X.
\end{equation}

Since we have assumed that  $X_n$ or $P_n$  is orthogonally (or unitarily) invariant and that they are independent, this implies that $U_n$ is a Haar-distributed orthogonal (or unitary) matrix that is also independent of    $(\la_1,\ldots, \la_n)$ (see the first paragraph of the proof of \cite[Th. 4.3.5]{hp00} for additional details). 

Recall that the largest eigenvalue $\la_1(X_n)\convas b$, while the smallest eigenvalue $\la_n(X_n)\convas a$. Let us now consider the eigenvalues of $\wtX$ which are out of $[\la_n(X_n), \la_1(X_n)]$. By Proposition \ref{9709.10h}-a) and an application of the identity in Proposition \ref{9709.10h}-b) these eigenvalues are precisely the numbers $z \notin [\la_n(X_n), \la_1(X_n)]$ \st the $r\times r$ matrix
\begin{equation}\label{eq:Mn z}
M_{n}(z):=I_r-[\theta_jG_{\mu_{i,j}^{(n)}}(z)]_{i,j=1}^r
\end{equation}
is singular.  Recall that in (\ref{eq:Mn z}), $G_{\mu_{i,j}^{(n)}}(z)$, for $i,j = 1, \ldots, r$ is the Cauchy transform of the random complex measure defined by  \be\label{1442010.1}
\mu_{i,j}^{(n)}=\sum_{m=1}^n\ovl{u_{k,i}^{(n)}}u_{k,j}^{(n)}\delta_{\la_k(X_n)},
\ee
where $u_{k,i}$ and $u_{k,j}$ are the $(k,i)-$th and $(k,j)$-th entries of the orthogonal (or unitary) matrix $U_n$ in \eqref{eq:the equivalent matrix 1} and $\lambda_{k}$ is the $k$-th largest eigenvalue of $X_n$ as in the first term in (\ref{eq:the equivalent matrix 1}).

By Proposition \ref{120709.8h35}, we have that
\begin{equation}
\mu_{i,j}^{(n)} \convas
\begin{cases}
\mu_{X} & \textrm{ for } i=j, \\
\delta_{0} & \textrm{otherwise.}
\end{cases}
\end{equation}

Thus we have
\begin{equation}\label{eq:Mnz uniform}
M_{n}(z) \overset{\textrm{a.s.}}{\longrightarrow} M_{G_{\mu_X}}(z):=\Diag(1-\theta_1G_{\mu_X}(z), \ldots, 1-\theta_rG_{\mu_X}(z)),
\end{equation}
and by Lemma \ref{140709.07h55},  this convergence is uniform on $\{z\in \C\ste \operatorname{d}(z, [a,b])\geq \eta >0\}$ for each $\eta>0$.

We now note that Hypotheses a), b) and c) of Lemma \ref{140709.19h26} are satisfied and follow from the definition of the Cauchy transform $G_{\mu_X}$. Hypothesis d) of Lemma \ref{140709.19h26} follows from the fact that $\wtX$ is Hermitian while hypothesis e) has been established in (\ref{eq:Mnz uniform}).

Let us recall that the eigenvalues of $\widetilde{X}_{n}$ which are out of $[\la_n(X_n), \la_1(X_n)]$ are precisely those values $z_n$ where the matrix $M_n(z_n)$ is singular. As a consequence, we are now in a position where Theorem \ref{140709.main} follows by invoking Lemma \ref{140709.19h26}.   Indeed, by Lemma \ref{140709.19h26}, if $$z_1>\cdots >z_p$$ denote the solutions of the equation $$\qquad\qquad\prod_{i=1}^r(1-\tta_iG_{\mu_X}(z))=0 \qquad\qquad\textrm{($z\in \R\bck[a,b]$),}$$ then their exists some sequences $(z_{n,1})$,\ldots, $(z_{n,p})$ converging respectively to $z_1, \ldots, z_p$ \st for any $\eps>0$ small enough, for   $n$ large enough,   the eigenvalues of $\wtX$ that are out of $[a-\eps, b+\eps]$  are exactly $z_{n,1}, \ldots, z_{n,p}$.
Moreover,  \eqref{8.10.10.19h55} and Lemma \ref{140709.19h26}-(ii) ensure that for $n$ large enough, these eigenvalues  have multiplicity one.

\subsubsection{Isolated eigenvalues in the setting where some $\theta_i$'s may coincide}\label{71009.15h19}
We now treat the case where the $\theta_i$'s are not supposed to be pairwise distinct.

 We denote, for $\theta\ne 0$, $$\rho_\theta=  \begin{cases} G_{\mu_X}^{-1}(1/\theta)&\textrm{ if $G_{\mu_X}(a^-)<1/\theta  <G_{\mu_X}(b^+)$,}\\ \\
b &\textrm{ if $1/\tta>G_{\mu_X}(b^+)$,}\\ \\
a&\textrm{ if $1/\tta<G_{\mu_X}(a^-)$.}
\end{cases}$$
We want to prove that for all $1\le i\le s$, $\la_i(\wtX)\convas \rho_{\tta_i}$   and that for all $0\le j< r-s$,  $\la_{n-j}(\wtX)\convas \rho_{\tta_{r-j}}$.

We shall treat only the case of  largest eigenvalues (the case of smallest ones can be treated in the same way).  So let us fix $1\le i\le s$ and $\eps>0$.

There is $\eta>0$ \st $|\rho_\theta-\rho_{\theta_{i}}|\leq \eps$
whenever $|\theta-\theta_{i}|\leq \eta$.
Consider pairwise distinct non zero real numbers  $\theta'_{1}>  \cdots>\theta'_{r}$  \st for all $j=1, \ldots, r$, $\theta_j$ and  $ \theta'_{j}$ have the same sign
 and
\bes\label{230709.06h57}\sum_{j=1}^r(\theta'_{j}-\theta_j)^2\leq\min(\eta^2,\eps^2).\ees
It implies that $|\rho_{\theta'_{i}}-\rho_{\theta_{i}}|\leq \eps$. 
With the notation in Section \ref{section18.09.09.10h35}, for  each $n$, we define \bes P'_{n}=U_P\Diag(\theta'_{1},\ldots, \theta'_{r}, 0, \ldots, 0)U_P.\ees
Note that by \cite[Cor. 6.3.8]{hj85}, we have, for all $n$, \bes\label{140709.19h58}\sum_{j=1}^n( \la_j(X_n+P'_{n}) -\la_j(X_n+P_n))^2\leq \Tr(P'_{n}-P_n)^2=\sum_{j=1}^r(\theta'_{j}-\theta_j)^2 \leq \eps^2 .\ees   Theorem \ref{140709.main} can applied to $X_n+P'_{n}$ (because the $\tta_1', \ldots, \tta_r'$ are pairwise distinct). It follows that almost surely, for $n$ large enough, \bes\label{23709.7h17}|\la_{i}(X_n+P'_n)-\rho_{\theta'_{i}}|\leq \eps.\ees By the triangular inequality,
almost surely, for $n$ large enough, \bes |\la_{i}(X_n+P_n)-\rho_{\theta_{i}}|\leq 3\eps,\ees
so that $\lambda_{i}(X_n + P_n) \convas \rho_{\theta_i}$.{\hfill$\square$}



\section{Proof of Theorem \ref{180709.13h39}} As in Section \ref{section18.09.09.10h35}, let \bes\label{140709.19h47lendemain}X_n=U_X\Diag(\la_1,\ldots, \la_n)U_X^{*}\, ,\qquad  P_n=U_P\Diag(\theta_1,\ldots, \theta_r, 0, \ldots, 0)U_P^{*}.\ees

The eigenvectors  of $X_n+P_n$, are precisely $U_X$ times the eigenvectors of $$\Diag(\la_1,\ldots, \la_n)+\underbrace{U_X^{*}U_P}_{= U_n}\Diag(\theta_1,\ldots, \theta_r, 0, \ldots, 0)U_P^{*}U_X.$$ Consequently, we have proved Theorem \ref{180709.13h39} by proving the result in the setting where $X_n= \Diag(\la_1,\ldots, \la_n)$ and $P_n= U_n\Diag(\theta_1,\ldots, \theta_r, 0, \ldots, 0)U_n^{*}$, where $U_n$ is a Haar-distributed orthogonal (or unitary) matrix.

As before, we denote  the entries of $U_n$ by $[u_{i,j}^{(n)}]_{i,j=1}^n$. Let the columns of $U_n$ be denoted by  $u_{1}, \ldots, u_n$ and the $n\times r$ matrix whose columns are respectively $u_1$, \ldots, $u_r$ by $U_{n,r}$. Note that the $u_i$'s, as $\widetilde{u}$, depend on $n$, hence should be denoted for example by $u_i^{(n)}$ and $\widetilde{u}^{(n)}$. The same, of course, is true for the $\la_i$'s and the $\widetilde{\la}_i$'s. To simplify our notation, we shall suppress the $n$ in the superscript.

  Let $r_0$ be the number of $i$'s \st $\theta_i=\theta_{i_0}$. Up to a reindexing of the $\theta_i$'s (which are then no longer decreasing - this fact does not affect our proof), one can suppose that ${{i_0}}=1$,  $\theta_1=\cdots=\theta_{r_0}$. This choice implies that, for each $n$, $\ker (\theta_1I_n-P_n)$ is    the linear span of the ${r_0}$ first columns $u_1$, \ldots, $u_{r_0}$ of $U_n$. By construction, these columns are orthonormal. Hence, we will have proved Theorem  \ref{180709.13h39} if we can prove that as  $n \longrightarrow \infty$,  \be\label{16709.18h28} \sum_{i=1}^{r_0}| \lan  u_{i}, \widetilde{u} \ran|^2 \convas \f{-1}{\theta_{i_0}^2G_{\mu_X}'(\rho)}=\ff{\theta_{i_0}^2\int \f{\ud \mu_{X}(t)}{(\rho-t)^2}} \ee and  \be\label{16709.18h28prime} \sum_{i={r_0}+1}^r| \lan  u_{i},\widetilde{u}\ran|^2 \convas 0. \ee

As before, for every $n$ and for all $z$ outside the spectrum of $X_n$,  define the $r\times r$ random matrix: $$M_{n}(z):=I_r-[\theta_jG_{\mu_{i,j}^{(n)}}(z)]_{i,j=1}^r,$$   where, for all $i,j=1, \ldots, r$, $\mu_{i,j}^{(n)}$ is the random complex measure defined by   \eqref{1442010.1}. 

In (\ref{eq:Mnz uniform}) we established that:
\be\label{271010.14h28}M_{n}(\cdot) \convas M_{G_{\mu_X}}(\cdot):=\Diag(1-\theta_1G_{\mu_X}(\cdot), \ldots, 1-\theta_rG_{\mu_X}(\cdot))\ee
uniformly on  $\{z\in \C\ste \operatorname{d}(z, [a,b])\geq \eta\}$ for all $\eta >0$.

We have established in Theorem \ref{140709.main} that because $\theta_{i_0} > 1/G_{\mu_X}(b^+) $,  $\wtl_{i_0} \convas \rho=G_{\mu_X}^{-1}(1/\theta_{i_0}) \notin [a,b]$ as $n\lto\infty$. It follows that:
\be\label{1442010.7h05}
M_{n}(z_n) \convas \Diag(\underbrace{0,\ldots, 0}_{{r_0} \textrm{ zeros}}, 1-\f{\theta_{{r_0}+1}}{\theta_{i_0}}, \ldots,  1-\f{\theta_{r}}{\theta_{i_0}}).\ee

Proposition \ref{9709.10h} a) states that for $n$ large enough so that \st $\widetilde{\lambda}_{i_0}$   is not an eigenvalue of $X_n$, the $r \times 1 $ vector $$U_{n,r}^*\widetilde{u} =\begin{bmatrix} \lan  u_{1}, \widetilde{u} \ran\\ \vdots \\ \lan  u_{r}, \widetilde{u} \ran \end{bmatrix}$$ is in the kernel of the $r\times r$ matrix $M_{n}(z_n)$ with  $||U_{n,r}^*\widetilde{u} ||_{2} \leq 1$.

Thus by \eqref{271010.14h28}, any limit point of $U_{n,r}^*\widetilde{u} $ is in the kernel of the matrix on the right hand side of \eqref{1442010.7h05}, i.e. has its $r-r_0$ last coordinates equal to zero.

Thus \eqref{16709.18h28prime} holds and we have proved Theorem \ref{180709.13h39}-b). We now establish \eqref{16709.18h28}.


By \eqref{130410.1flo-zurich},  one has that for all $n$, the eigenvector $\widetilde{u}$ of $\wtX$ associated with the eigenvalue $\wtl_{i_0}$ can be expressed as: \beqy \nonumber \widetilde{u}&=&(\wtl_{i_0} I_n-X_n)^{-1}U_{n,r}\Diag(\theta_1, \ldots, \theta_r)U_{n,r}^*\widetilde{u}\\ \nonumber
&=& (\wtl_{i_0} I_n-X_n)^{-1}\sum_{j=1}^r\theta_j\lan u_j,\widetilde{u}\ran u_j,\\\nonumber
&=& \underbrace{(\wtl_{i_0} I_n-X_n)^{-1}\sum_{j=1}^{r_0}\theta_j\lan u_j,
\widetilde{u}\ran u_j}_{\textrm{denoted by }\widetilde{u}'}+ \underbrace{(\wtl_{i_0}I_n-X_n)^{-1}\sum_{j={r_1}+1}^r\theta_j\lan u_j,\widetilde{u}\ran u_j}_{\textrm{denoted by }\widetilde{u}''}.
\eeqy
As $\wtl_{i_0} \convas \rho \notin [a,b]$,  the sequence $ (\wtl_{i_0} I_n-X_n)^{-1}$ is bounded in operator norm so that by  \eqref{16709.18h28prime}, $\|\widetilde{u}''\|\convas 0$.  Since $\|\widetilde{u}\|=1$, this implies that  $\|\widetilde{u} '\|\convas 1$.

Since we assumed that $\tta_{i_0}=\theta_1=\cdots=\theta_{r_0}$, we must have that:
\beqy  \label{17709.14hprime}
 \|\widetilde{u}'\|^2
 &=&
\theta_{i_0}^2 \sum_{i,j=1}^{r_1}\ovl{\lan u_i, \widetilde{u}\ran} \lan u_{j}, \widetilde{u}\ran  \underbrace{u_i^{*}(z_n I_n-X_n)^{-2}u_{j}}_{=\int\f{1}{(z_n-t)^2}\ud\mu^{(n)}_{i,j}(t)}.
\eeqy

By Proposition \ref{120709.8h35}, we have that for all $i\neq j$, $\mu_{i,j}^{(n)} \convas \delta_{0}$ while for all $i$, $\mu_{i,i}^{(n)} \convas \mu_{X}$. Thus, 
since we have that $z_n \convas \rho \notin[a,b]$, we have that for all $i,j=1, \ldots, r_0$, $$ \int\f{\ud\mu^{(n)}_{i,j}(t)}{(z_n-t)^2}\convas \one_{i=j} \int\f{\ud\mu (t)}{(\rho-t)^2}.$$

Combining  the relationship in \eqref{17709.14hprime} with the fact that $\|\widetilde{u}'\|\convas 1$,  yields \eqref{16709.18h28} and we have proved Theorem \ref{180709.13h39}-a).{\hfill$\square$}

\section{Proof of Theorem \ref{200709.13h5}} Let us assume that $\theta > 0$. The proof supplied below can be easily ported to the setting where $ \theta < 0$.

We first note that the $G_{\mu_X}'(b^+)=-\infty$, implies that $\int\f{\ud\mu_X(t)}{(b-t)^2}= - G_{\mu_X}(b^{+}) = +\infty$.
We adopt the strategy employed in proving Theorem \ref{180709.13h39}, and note that it suffices to prove the result in the setting where $X_n= \Diag(\la_1,\ldots, \la_n)$ and $P_n=\theta uu*$, where $u$ is a $n \times 1$ column vector uniformly distributed on the unit real (or complex) $n$-sphere.

We denote the coordinates of $u$ by $u^{(n)}_1, \ldots, u^{(n)}_n$ and define, for each $n$, the random \pro measure $$\mu_{X}^{(n)}=\sum_{k=1}^n|u^{(n)}_k|^2\delta_{\la_k}.$$

The  $r = 1$ setting of Proposition \ref{9709.10h}-b) states that the eigenvalues of $X_n+P_n$ which are not eigenvalue of $X_n$    are the solutions of $$G_{\mu_{X}^{(n)}}(z)=\ff{\theta}.$$

Since $G_{\mu_{X}^{(n)}}(z)$ decreases from $+\infty$ to $0$ for increasing values of $z \in (\la_1, +\infty)$, we have that $\la_1(X_n+P_n) =: \widetilde{\lambda}_{1} >\la_1$. Reproducing the arguments leading to \eqref{250909.17h05} in Section \eqref{250909.17h06}, yields the relationship:
\begin{equation}\label{eq:phase zero ev}
|\lan \widetilde{u}, \ker(\theta I_n-P_n)\ran|^2=\ff{\theta^2\int  \f{\ud \mu^{(n)}(t)}{(\widetilde{\lambda}_1-t)^2}}.
\end{equation}
Thus, proving that:$${\int  \f{\ud \mu^{(n)}(t)}{(\widetilde{\lambda}_1-t)^2}}\convas+\infty \quad \textrm{ as } n \longrightarrow \infty$$ will yield the desired result.

By hypothesis, we have that $$\ff{n}\sum_{i=1}^n\delta_{\la_i}\convas \mu_X.$$ By Theorem \ref{140709.main}, we have that $\widetilde{\lambda}_1\convas b$ so that $$\ff{n}\sum_{k=1}^n\delta_{\la_k +b-\widetilde{\lambda}_1}\convas \mu_X.$$
Hence, by \eqref{1179.22h34a}, $$ \tilde{\mu}^{(n)}:=\sum_{k=1}^n|u^{(n)}_k|^2\delta_{\la_k+b-\widetilde{\lambda}_1}\convas \mu_X.$$
This implies that almost surely,
$$\liminf_{n\to \infty}{\int  \f{\ud \mu^{(n)}(t)}{(\widetilde{\lambda}_1-t)^2}}=\liminf_{n\to \infty}\int \f{\ud \tilde{\mu}^{(n)}(t)}{(b-t)^2}\ge\int \f{\ud {\mu_X}(t)}{(b-t)^2} =+\infty,$$
so that by (\ref{eq:phase zero ev}), $\lan \widetilde{u}, \ker(\theta I_n-P_n)\ran|^2 \convas 0$ thereby proving Theorem \ref{200709.13h5}.{\hfill$\square$}
\\ \\

We omit the details of the proofs of Theorems \ref{140709.main.multiplicative}--\ref{200709.13h5.multiplicative}, since these are straightforward adaptations of the proofs of  Theorems \ref{140709.main}-Theorem \ref{200709.13h5} that can obtained by following the prescription in Proposition \ref{9709.10h}-c).

%% file: convergencefacts_fin_v6.tex
\section{Appendix: convergence of weighted spectral measures}\label{120709.8h34}

\subsection{A few facts about the weak convergence of complex measures}
Recall that a sequence $(\mu_n)$ of complex measures on  $\R$ is said to {\it converge weakly} to a complex measure $\mu$ on $\R$ if, for any continuous bounded function $f$ on $\R$, \be\label{10709.10h542}\int f(t)\ud \mu_n(t)\ninf \int  f(t)\ud \mu(t).\ee
We now establish a lemma on the weak convergence of complex measures that will be useful in proving Proposition \ref{120709.8h35}. We note that the counterpart of this lemma for \pro measures is well known. We did not find any reference in standard literature to the ``complex measures version'' stated next, so we provide a short proof.

Recall that a sequence $(\mu_n)$ of complex measures on $\R$ is said to be {\it tight} if $$\lim_{R\to+\infty}\sup_n |\mu_n|(\{t\in \R\ste |t|\geq R\})=0.$$
\begin{lem}\label{140709.07h55}
Let $D$ be a dense subset of  the set of continuous functions on $\R$ tending to zero at infinity endowed with the topology of the uniform convergence. Consider a tight sequence $(\mu_n)$ of complex measures on $\R$ \st   $(|\mu_n|(\R))$ is bounded and  \eqref{10709.10h542} holds for any  function $f$ in $D$. Then  $(\mu_n)$ converges weakly to $\mu$. Moreover, the convergence of \eqref{10709.10h542} is uniform on any set of uniformly bounded and uniformly Lipschitz functions.
\end{lem}


\begin{pr} Firstly, note that using the boundedness of $(|\mu_n|(\R))$, one can easily   extend \eqref{10709.10h542} to any continuous function tending to zero at infinity.
It follows that for any continuous bounded function $f$ and any  continuous function $g$ tending to zero at infinity, we have
$$\limsup_{n\to\infty} \lf|\int f\ud (\mu-\mu_n)\ri|\leq \sup_{n}  \int |f (1-g)|\ud (|\mu|+|\mu_n|),$$ which can be made arbitrarily small by appropriately choosing $g$. The tightness hypothesis ensures that such a $g$ can always be found. This proves that $(\mu_n)$ converges weakly to $\mu$. The uniform convergence follows from a straightforward application of Ascoli's Theorem.\end{pr}

\subsection{Convergence of weighted spectral measures}

\begin{lem}\label{14.10.2010.16h52}For each $n$, let $u^{(n)}=(u_{1},\ldots, u_{n})$,  $v^{(n)}=(v_1,\ldots, v_n)$ be the first two rows of a Haar distributed random orthogonal (or unitary) matrix. Let  $x^{(n)}=(x_1, \ldots,x_n)$ be a set of real-valued random numbers\footnote{\label{8.11.2010.1}Note that for each $k$, $u_k,v_k$ and $x_k$ obviously depend on $n$. We could have represented this dependence in the subscript as $u_{n,k}$ or $u_k^{(n)}$ but we choose to suppress the $n$ for notational brevity.} independent of $(u^{(n)},  v^{(n)})$. We suppose that $\sup_{n,k}|x_k|<\infty$ almost surely.

a)  Then as $n \longrightarrow \infty$, \bes\label{210709.12h31bis}  \ovl{u}_1 v_1 x_1 + \ovl{u}_2 v_2 x_2 +\cdots+ \ovl{u}_n v_n x_n \convas 0.\ees

b)  Suppose   that $\ff{n}(x_1 +x_2 +\cdots+x_n )$ converges almost surely to a deterministic limit $l$. Then \bes\label{210709.12h32}  |u_1 |^2x_1 +|u_2 |^2x_2 +\cdots+|u_n |^2x_n \convas l.\ees
\end{lem}

\begin{pr}
We prove a) by conditioning on the values of the $x $'s. As a consequence, we can suppose the $x^{(n)} $'s to be non-random. Let us define the random variable $Z_n=\ovl{u}_1 v_{1} x_1 +\cdots+ \ovl{u}_n {v_n} x_n $. To establish a), it suffices to prove that $\E[Z_n^4]=O(n^{-2})$. We have
 $$\E[Z_n^4]=
 \sum_{i,j,k,l=1}^n
 x_i x_j x_k x_l  \E[u _i u _ju _k u _l v _iv _jv _k v _l   ].$$ Since $(u _{1},\ldots, u _{n})$ and $(v _1,\ldots, v _n)$ are the first two rows of a Haar distributed random orthogonal (or unitary) matrix, we have that $\E[u _i u _ju _k u _l v _iv _jv _k v _l   ]=0$ whenever one of the indices $i,j,k,l$, is different from the others \cite{cs06}. It follows that  \beq \E[Z_n^4]&\leq& 3 \sum_{i,j}  x_i^{2} x_j^{2}\E \left[u _i^{2}  u _j ^{2}  v _i  ^{2}  v _j  ^{2}   \right]\\ &\leq& 3 \sum_{i,j} x_i^{2}x_j^{2}\underbrace{ \lf(\E\left[ u _i ^{8}\right]\E\left[  u _j ^{8} \right]\E\left[  v _i ^{8}\right]\E\left[  v _j ^{8}   \right] \ri)^\ff{4}}_{=\E\left[ u _1 ^{8}\right]}.
\eeq
Since, by \cite{cs06}, $\E[  u _1 ^{8} ]=O(n^{-4})$ and we have assumed that $\sup_{n,k}|x_k |<\infty$, $ \ovl{u}_1 {v_1} x_1 +\cdots+ \ovl{u}_n {v_n} x_n \convas 0$.

To prove b), we employ a different strategy. In the setting where $(u _{1},\ldots, u _{n})$ is the first row of a Haar distributed orthogonal matrix, the result follows from the application of a well-known concentration of measure result \cite[Th. 2.3 and Prop. 1.8]{ledoux-amsbook} which states that there are positive constants $C,c$ \st for $n$ large enough, for any $1$-Lipschitz function $f_n$ on the unit sphere of $\R^n$, for all $\eps >0$,
\bes\label{127092} \Pro\lf\{|f_n(u _{1},\ldots, u _{n})-\E[f_n(u _{1},\ldots, u _{n})]|\geq\eps \ri\}\leq Ce^{-cn\eps^2}. \ees This implies that if $\E[f_n(u _{1},\ldots, u _{n})] $ converges, as $n \longrightarrow \infty$, to a finite limit, then $f_n(u _{1},\ldots, u _{n})$ converges almost surely to the same limit. In the  complex  setting, we note that a uniformly distributed random vector on the unity sphere of $\C^n$ is a uniformly distributed random vector on the unit sphere of $\R^{2n}$ so that we have proved the results in the unitary setting as well.\end{pr}


\begin{propo}\label{120709.8h35} Let, for each $n$, $u^{(n)} =(u _1,\ldots, u _n)$,  $v^{(n)}= (v _1,\ldots, v _n)$ be the two first columns of a uniform random orthogonal (resp. unitary) matrix.  Let also $\la^{(n)} =(\la _1, \ldots,\la _n)$ be a random family of real numbers\footnote{As in Lemma \ref{14.10.2010.16h52}, we have suppressed the index $n$ for notational brevity.} independent of  $(u^{(n)},  v^{(n)})$ \st almost surely, $\sup_{n,k}|\la_k|<\infty$.  We suppose that    there exists a deterministic \pro measure $\mu$ on $\R$  \st  almost surely, as $n \longrightarrow 0$,  $\ff{n}\sum_{k=1}^n\delta_{\la _k}$ converges weakly to $\mu$.


 Then as $n \longrightarrow \infty$,
\beqy\label{1179.22h34a}&\mu_{U^{(n)}}:= \sum_{k=1}^n|{u_k }|^2\delta_{\la_k} \textrm{ converges almost surely  weakly to $\mu$,}&\\
&\label{1179.22h34}\mu_{U^{(n)}, V^{(n)}}:=\sum_{k=1}^n\ovl{u}_k v_k \delta_{\la_k} \textrm{ converges almost surely weakly to $0$.}&
\eeqy
\end{propo}

\begin{pr}
We   use  Lemma \ref{140709.07h55}. Note first that   almost surely, since $\sup_{n,k}|\la_k|<\infty$,  both  sequences   are tight. Moreover, we have $$|\mu_{U^{(n)}, V^{(n)}}|=\sum_{k=1}^n|u_kv_k|\delta_{\la_k},$$ thus, by the Cauchy-Schwartz inequality,   $|\mu_{U^{(n)}, V^{(n)}}|(\R)\leq 1$. 
 The set of continuous functions on the real line tending to zero at infinity  admits a countable dense subset, so it suffices to prove that for any fixed such function $f$, the convergences of \eqref{1179.22h34a} and 
 \eqref{1179.22h34}
 hold almost surely when applied to $f$. This follows easily from Lemma \ref{14.10.2010.16h52}. \end{pr}